\documentclass[a4paper,10pt]{article}
\usepackage{amsmath}
\usepackage{amssymb}
\usepackage{amsthm}
\usepackage{amscd}
\begin{document}
\theoremstyle{plain}
\newtheorem{thm}{Theorem}[section]
\newtheorem{theorem}[thm]{Theorem}
\newtheorem{lemma}[thm]{Lemma}
\newtheorem{corollary}[thm]{Corollary}
\newtheorem{proposition}[thm]{Proposition}
\newtheorem{addendum}[thm]{Addendum}
\newtheorem{variant}[thm]{Variant}
\theoremstyle{definition}
\newtheorem{notations}[thm]{Notations}
\newtheorem{question}[thm]{Question}
\newtheorem{problem}[thm]{Problem}
\newtheorem{remark}[thm]{Remark}
\newtheorem{remarks}[thm]{Remarks}
\newtheorem{definition}[thm]{Definition}
\newtheorem{claim}[thm]{Claim}
\newtheorem{assumption}[thm]{Assumption}
\newtheorem{assumptions}[thm]{Assumptions}
\newtheorem{properties}[thm]{Properties}
\newtheorem{example}[thm]{Example}

\numberwithin{equation}{section}

\newcommand{\sA}{{\mathcal A}}
\newcommand{\sB}{{\mathcal B}}
\newcommand{\sC}{{\mathcal C}}
\newcommand{\sD}{{\mathcal D}}
\newcommand{\sE}{{\mathcal E}}
\newcommand{\sF}{{\mathcal F}}
\newcommand{\sG}{{\mathcal G}}
\newcommand{\sH}{{\mathcal H}}
\newcommand{\sI}{{\mathcal I}}
\newcommand{\sJ}{{\mathcal J}}
\newcommand{\sK}{{\mathcal K}}
\newcommand{\sL}{{\mathcal L}}
\newcommand{\sM}{{\mathcal M}}
\newcommand{\sN}{{\mathcal N}}
\newcommand{\sO}{{\mathcal O}}
\newcommand{\sP}{{\mathcal P}}
\newcommand{\sQ}{{\mathcal Q}}
\newcommand{\sR}{{\mathcal R}}
\newcommand{\sS}{{\mathcal S}}
\newcommand{\sT}{{\mathcal T}}
\newcommand{\sU}{{\mathcal U}}
\newcommand{\sV}{{\mathcal V}}
\newcommand{\sW}{{\mathcal W}}
\newcommand{\sX}{{\mathcal X}}
\newcommand{\sY}{{\mathcal Y}}
\newcommand{\sZ}{{\mathcal Z}}
\newcommand{\bbA}{{\mathbb A}}
\newcommand{\bbB}{{\mathbb B}}
\newcommand{\bbC}{{\mathbb C}}
\newcommand{\bbD}{{\mathbb D}}
\newcommand{\bbE}{{\mathbb E}}
\newcommand{\bbF}{{\mathbb F}}
\newcommand{\bbG}{{\mathbb G}}
\newcommand{\bbH}{{\mathbb H}}
\newcommand{\bbI}{{\mathbb I}}
\newcommand{\bbJ}{{\mathbb J}}
\newcommand{\bbM}{{\mathbb M}}
\newcommand{\bbN}{{\mathbb N}}
\newcommand{\bbQ}{{\mathbb Q}}
\newcommand{\bbR}{{\mathbb R}}
\newcommand{\bbT}{{\mathbb T}}
\newcommand{\bbU}{{\mathbb U}}
\newcommand{\bbV}{{\mathbb V}}
\newcommand{\bbW}{{\mathbb W}}
\newcommand{\bbX}{{\mathbb X}}
\newcommand{\bbY}{{\mathbb Y}}
\newcommand{\bbZ}{{\mathbb Z}}
\newcommand{\id}{{\rm id}}
\newcommand{\rank}{{\rm rank}}

\title{\textbf{Effective behavior of multiple linear systems}
\thanks{The author would like to thank the hospitality and
financial support of The Institute of Mathematical Sciences in
The Chinese University of Hong Kong during this research. This work
is also supported by 973 Science Foundation, the Doctoral Program
Foundation of EMC and the Foundation of Shanghai for Priority
Academic Discipline.}}
\author{{Sheng-Li Tan}
\\[8pt]
Department of Mathematics \\
East China Normal University\\
Shanghai 200062, P.~R. of China\\
sltan@math.ecnu.edu.cn}

\def\d{\partial}
\def\c_t{\chi_{\text{top}}}
\def\BD{\boldsymbol D}
\def\BDelta{\boldsymbol\Delta}
\def\wt{\widetilde}
\def\wh{\widehat}
\def\la{\leftarrow}
\def\ep{\varepsilon}
\def\ra{\rightarrow}
\def\lra{\longrightarrow}
\def\G{\Gamma}
\def\O{\mathcal O}
\def\I{\mathcal I}
\def\L{\mathcal L}
\def\M{\mathcal M}
\def\E{\mathcal E}
\def\F{\mathcal F}
\def\H{\mathcal H}
\def\l{\ell}
\def\lra{\longrightarrow}
\def\ra{\rightarrow}
\def\Pic{\text{\rm{Pic}}}
\def\Spec{\text{\rm{Spec\,}}}
\def\grad{\text{\rm{grad\,}}}
\def\syz{\text{\rm{syz\,}}}
\def\tr{\text{\rm{tr\,}}}
\def\div{\text{\rm{div}}}
\def\Div{\text{\rm{Div}}}
\def\parlistitem#1{\par\hskip0.7cm\llap{(#1)}\ }
\def\ratmap{\,-\,-\to\,}
\def\phi{\varphi}
\def\({$($}
\def\){$)$}
\def\disc{\text{\rm disc\,}}

\font\eightrm=cmr10 at 8pt
\font\bigrm=cmr10 scaled1440
\font\midrm=cmr10 scaled1200
\font\ninerm=cmr10 at 9pt
\font\ninebf=cmbx10 at 9pt
\font\bbf=msbm10
\font\bigbf=cmbx10 scaled1440
\font\Bigbf=cmbx10 scaled2073
\font\midbf=cmbx10 scaled1200
\font\eightit=cmti8
\font\nineit=cmti9
\renewcommand{\labelenumi}{$($\arabic{enumi}$)$}
\renewcommand{\labelenumii}{(\Alph{enumi})}
\renewcommand{\atop}[2]{\genfrac{}{}{0pt}{}{#1}{#2}}
\date{Dedicated to Prof. Yum-Tong Siu on his 60th birthday}

\maketitle

\section{Introduction}

It is a fundamental problem in algebraic geometry to understand the behavior of a multiple
linear system $|nD|$ on a projective complex manifold $X$ for large $n$. For example, the
well-known Riemann-Roch problem is to compute the function
$$
n\longmapsto h^0(\O_X(nD)):=\dim_{\mathbb C} H^0(X,\O_X(nD)).
$$
In the introduction to his collected works \cite{Zar2}, Zariski
cited the Riemann-Roch problem as one of the four ``difficult
unsolved questions concerning projective varieties (even algebraic
surfaces)''. The other natural problems about $|nD|$ are to find
the fixed part and base points (see \cite{Zar1}), the very
ampleness, the properties of the associated rational map and its
image variety, the finite generation of the ring of sections,
$\cdots$.

For a genus $g$ curve $X$, Riemann-Roch theorem gives good and {\it effective} solutions to these problems.
\begin{enumerate}
\item[$\bullet$] Assume that $\deg D>0$. If $n\geq \frac{2g-1}{\deg D}$, \ then $h^1(\O_X(nD))=0$.
So, in general,
\begin{equation*}
h^0(nD)=\begin{cases} n\deg D+1-g, &\textrm{ if } n\deg D>2g-2,\\
\textrm{a periodic function of } n, &\textrm{ if } \deg D=0,\\
0, &\textrm{ if } \deg D<0.
\end{cases}
\end{equation*}
\item[$\bullet$] If $n\geq \frac{2g}{\deg D}$, \ then $|nD|$ is base point free.
\item[$\bullet$] If $n\geq \frac{2g+1}{\deg D}$, \ then $|nD|$ is very ample.
\end{enumerate}

When $X$ is a surface, the Riemann-Roch problem is also equivalent to the computation of
$h^1(\O_X(nD))$. This problem was studied first by the Italian geometers in the 19th century.
Castelnuovo \cite{Cas} proved that if $|D|$ is a base point
free linear system of dimension $\geq 2$, then there is a constant $s$ such that
\begin{equation*}
h^0(\O_X(nD))=\chi(\O_X(nD))+s
\end{equation*}
for $n$ sufficiently large, i.e., $h^1(\O_X(nD))$ is a constant.

In \cite{Zar1}, Zariski established the fundamental theory on the
behavior of an arbitrary multiple linear system $|nD|$ on an
algebraic surface (see the next section for the details.) By using
{\it Zariski decomposition}, he showed that one only needs to know
the behavior of $|nA+T|$, where $A$ is a nef divisor and $T$ is a
fixed divisor (see Theorem \ref{THM}). Zariski proved the
boundedness of $h^1(\O_X(nA+T))$, the fixed part $B_n$ and the
isolated base points of $|nA+T|$ when $n$ is sufficiently large.
An important conjecture on the periodicity was proved later by
Cutkosky and Srinivas \cite{CuS1} in 1993. However, all of these
results are ineffective on $n$.

In the language of Beltrametti and Sommese \cite{BeS1}, these
problems are about the $k\,$-very ampleness.

\begin{definition} (1) Let $k$ be a nonnegative integer. A divisor $D$
(or the linear system $|D|$) on $X$ is called
$k$-very ample if any $k+1$ points (not necessarily distinct) give $k+1$ independent
conditions on $|D|$. Precisely, for any zero dimensional subscheme $\Delta\subset X$ with
$\deg\Delta:=h^0(\mathcal O_\Delta)\leq k+1$,
\begin{equation*}
h^0(I_\Delta(D))=h^0(\O_X(D))-\deg\Delta,
\end{equation*}
where $I_\Delta$ is the ideal sheaf of $\Delta$.

(2) $D$ is called $(-1)$-very ample if $H^1(\O_X(D))=0$.
\end{definition}

The Riemann-Roch problem is about $(-1)$-very ampleness;
``$0$-very ample'' is equivalent to ``base point free'';
``$1$-very ample'' is just ``very ample''.

If $X$ is a curve of genus $g$, $k\geq -1$ and $D$ is a divisor
satisfying
$$\deg D\geq k+2g,$$
 then $D$ is $k$-very ample. In particular, if $\deg D>0$, then $|nD|$
is $k$-very ample provided
$$
n\geq \dfrac{k+2g}{\deg D}.
$$

In recent years, the effective version of some important theorems
attracted much attention. For example, Fujita's conjecture and the
effective Matsusaka's big theorem (see, for example, \cite{Siu1,
Siu2, Siu96, Siu3, Dem, Dem2, EL, Fuj, Kaw, Hei} ...).
 They are about the 0- and 1-very ampleness of the adjoint
linear system $|nH+K_X|$ and $|nH|$ for an ample divisor $H$, here
$K_X$ is the canonical divisor of $X$. I would like to mention the
latest bounds of Angehrn-Siu \cite{Siu1} and Siu \cite{Siu3} for a $d$-dimensional
complex manifold: $|K_X+nH|$ is base point free if
$$
n\geq\dfrac12(d^2+d+2).
$$
$|nH|$ is very ample if
$$
n\geq
\dfrac{\left(2^{3d-1}5d\right)^{4^{d-1}}\left(3(3d-2)^dH^d+K_XH^{d-1}\right)^{4^{d-1}3d}
}
{\left(6(3d-2)^d-2d-2\right)^{4^{d-1}d-\frac23}\left(H^d\right)^{4^{d-1}3(d-1)}}.
$$

 If $X$ is a surface, then there are also nice solutions: $|nH+K_X|$ is $(n-3)$-very
ample (Reider \cite{Rei}, Beltrametti and Sommese \cite{BeS1}). $|nH|$
is very ample when
\begin{equation}\label{equ1.1}
n>\dfrac12\left[\dfrac{(H(K_X+4H)+1)^2}{H^2}+3\right],
\end{equation}
(Fernandez del Busto \cite{Fer}). This bound is improved by Beltrametti and Sommese \cite{BeS2}
\begin{equation}\label{equ1.2}
n>\dfrac12\left(\dfrac{(H(K_X+2H)+1)^2}{H^2}+7\right).
\end{equation}
But the optimal effective Matsusaka's theorem for a surface is still open
(Ein \cite{Ein}, Open Problem 4).

It is also of great interest to find the effective behavior of a multiple linear system $|nD|$.
The purpose of this note is to give effective version of some well-known
theorems on multiple linear systems due to Zariski \cite{Zar1}, Castelnuovo \cite{Cas},
 Artin \cite{Artin1, Artin2}, Benveniste \cite{Ben}, Cutkosky and
Srinivas \cite{CuS1, CuS2}. We also try to find the effective
behavior of the rational map defined by $|nD|$.

For two divisors $A$ and $T$, let
\begin{align*}
\frak M(A,T)&=\dfrac{((K_X-T)A+2)^2}{4A^2}-\dfrac{(K_X-T)^2}{4},\\
\frak m(A,T)&=\min\{\, n \in \mathbb Z \ | \ n>\frak M(A,T)\}.
\end{align*}
Now we state our main result.

\begin{theorem}\label{Main Theorem} Let $A$ be a nef and big divisor on an
algebraic surface $X$, let
$T$ be any fixed divisor, and let $k$ be a nonnegative integer.
Assume that either $n>k+\frak M(A,T)$,  or $n\geq \frak M(A,T)$
when $k=0$ and $T\sim K_X+\lambda A$ for some $\lambda\in \mathbb
Q$. Suppose $|nA+T|$ is not $(k-1)$-very ample, i.e., there is a
zero dimensional subscheme $\Delta$ on $X$ with minimal degree
$\deg\Delta\leq k$ such that it does not give independent
conditions on $|nA+T|$. Then there is an effective divisor $D\neq
0$ containing $\Delta$ such that
\begin{equation}\label{equ1.3}
\begin{cases}
TD-D^2-K_XD\leq k, &\cr DA=0.
\end{cases}
\end{equation}
\end{theorem}

The effective version of the theorems in \cite{Zar1, Artin1,
Artin2, Ben, Cas, CuS1, CuS2} are obtained by direct applications
of this theorem for various $T$. For example, if $A=H$ is ample,
then $(\ref{equ1.3})$ has no solution $D\neq 0$. Thus we get an
effective version of Serre's theorem.
\begin{corollary}\label{Matsusaka}
If $H$ is ample and $n>k+\frak M(H,T)$ for some integer $k\geq 0$,
then $|nH+T|$ is $(k-1)$-very ample. Equivalently, $|nH+T|$ is
$(n-\frak m-1)$-very ample when $n\geq \frak m$.
\end{corollary}

The bound in this corollary is optimal in many cases. If $T=K_X$, then
$$
\frak M(A,K_X)=\dfrac{1}{A^2}.
$$
So $\frak m(A,K_X)=2$ (or $1$ if $A^2>1$). Thus the corollary
implies also that $|nH+K_X|$ is $(n-3)$-very ample (or
$(n-2)$-very ample if $H^2>1$).

If $T=0$, then
$$
\frak M(A,0)=\dfrac{(K_XA+2)^2}{4A^2}-\dfrac{K_X^2}4,
$$
and this corollary for $k=2$ is an effective version of Matsusaka
Big Theorem. Our bound is better than (\ref{equ1.1}) and
(\ref{equ1.2}). We will present an example to show that this bound
is the best possible.

In general, we set
\begin{equation*}
\tau(A,T)=\begin{cases}
 \underset{D}{\min}\,\{\,TD-K_XD-D^2\,\}, & \textrm{ if }A \textrm{ is not
 ample,}\cr
 +\infty, &\textrm{ if } A \textrm{ is
ample,}
\end{cases}
\end{equation*}
where $D$ runs over all effective divisors $D\neq 0$ such that
$DA=0$. $\tau$ is well defined (Lemma \ref{Lemma 2.7.4}). Then we
have
\begin{corollary} Assume that $\tau=\tau(A,T)\geq 1$ and $n\geq \frak m=\frak m(A,T)$.
Then $|nA+T|$ is $$\min\{\,\tau-2,\ n-\frak m-1\,\}$$ very ample.
\end{corollary}

Some well-known conditions on linear systems are those satisfying
$\tau\geq 1$ (see \S\,2). For example, $\tau(A,K_X)=-D^2\geq 1$
(Fujita's condition). $\tau(A,0)=2$ if and only if $p_a(D)\leq 0$
for any $D$ (Artin's condition). Laufer-Ramanujan's condition is
that $TD\geq K_XD$ for any $D$ (reduced and irreducible), which
implies also that $\tau\geq 1$.

As a consequence, the behavior of $|nA|$ is controlled by the
curves $C_i$ orthogonal to $A$, namely $AC_i=0$. If Artin's
condition is satisfied, then the behavior of $|nA|$ is quite
similar to that of the canonical multiple linear system $|nK_X|$
of a minimal surface of general type.

I would like to thank the referee for the valuable suggestions for
the correction of the original version.

\section{Zariski's results and generalizations}

In this section, we recall Zariski's fundamental results and their
generalizations. In our language, these results are essentially
about $(-1)$- and $0$-very ampleness.

 Let $X$ be a smooth projective complex surface, $K_X$ be its canonical
divisor and $D$ be any divisor on $X$.
\begin{definition} $D$ is called {\it nef } (numerically effective) if for any curve $C$ on
$X$, $DC\geq 0$. $D$ is called {\it big} if $D^2>0$. $D$ is called
{\it pseudo-effective} if for any ample divisor $H$, $DH\geq 0$.
\end{definition}

\begin{thm}\label{Zariski decomposition} {\bf (Zariski decomposition {\rm\cite{Zar1}}).}
 Let $D$ be a
pseudo-effective divisor on $X$. There exist uniquely $\mathbb Q$-divisors $A$ and
$F$ on $X$, such that $D=A+F$ satisfying the following conditions:
\begin{enumerate}
\item[\rm(1)] $F=0$ or the intersection matrix of the irreducible components of $F$
is negative definite;
\item[\rm(2)] $A$ is nef and $F$ is effective;
\item[\rm(3)] each irreducible component $C$ of $F$ satisfies $AC=0$.
\end{enumerate}
\end{thm}

The decomposition is called {\it Zariski decomposition}. The following basic theorem
has been used to reduce the general case $|nD|$ to the case
$|nA+T|$, where $T$ is any fixed divisor on $X$.

\begin{thm}\label{THM} {\rm({\bf Zariski} \cite{Zar1}).} As in Theorem \ref{Zariski decomposition}, $D=A+F$ is the Zariski decomposition.
\begin{enumerate}
\item[\rm(1)] $\kappa(D)=2$ if and only if $A^2>0$.
\item[\rm(2)] If $D$ is effective, then for all $n\geq 0$,
$$h^0(\O_X(nD))=h^0(\O_X([nA]))$$
\item[\rm(3)] If $sA$ is an integral divisor, and $n=as+b$ with $0\leq b<s$, then
$$h^0(\O_X(nD))=h^0(\O_X(asA+bD)).$$
\item[\rm(4)] As in {\rm(3)}, if $\kappa(D)\geq 0$, then
$$
h^1(\O_X(nD))=h^1(\O_X(asA+bD))-\dfrac{F^2}{2}{n^2}+\dfrac{FK_X}{2}n
-\left( b\dfrac{FK_X}{2}-b^2\dfrac{F^2}{2} \right).
$$
\end{enumerate}
\end{thm}

\begin{thm} {\rm({\bf Zariski} \cite{Zar1}).} As in Theorem \ref{THM}, assume that $T$ is any divisor on $X$.
\begin{enumerate}
\item[\rm(1)] $h^1(\O_X(nA+T))$ is bounded.
\item[\rm(2)] Let $D$ be effective and let $B_n$ be the fixed part of $|nD|$.
Then
$$B_n=\wt B_n+nF,$$
where $\wt B_n$ is a bounded
\(rational\,\) effective divisor.
\item[\rm(3)] If $|D|$ has no fixed part and $n\gg0$, then $|nD|$ has
no base point and $h^1(\O_X(nA+T))$ is a constant.
\end{enumerate}
\end{thm}

Shafarevich \cite{Sha} gave a new proof of the base point freeness
in (3).

\begin{thm} {\bf (Zariski {\rm\cite{Zar1}}, Cutkosky-Srinivas {\rm\cite{CuS1, CuS2}}).}
$h^1(\O_X(nA+T))$ and $\wt B_n$ are periodic when $n\gg0$.
\end{thm}

This theorem has been proved by Zariski \cite{Zar1} for the case $A^2=0$ and by
Cutkosky and Srinivas \cite{CuS1, CuS2} for the case $A^2>0$.

For a fixed $D$, we let
$$
R_m=R_m[D]=H^0(X,\mathcal O_X(mD)), \hskip0.3cm
R[D]=\oplus_{m=0}^\infty R_m[D].
$$
$R[D]$ is a graded ring.

Zariski gave in \cite{Zar1} a criterion for $R[D]$ to be finitely
generated.

\begin{thm}\label{generation} {\bf (Zariski {\rm\cite{Zar1}})} $R[D]$ is finitely
generated if and only if $\kappa(D)\leq 1$,
or $\kappa(D)=2$ and some multiple $|h(D-F)|$ has no fixed part.
\end{thm}

\begin{definition}\label{rat} (1) A curve $C=C_1+\cdots+C_r$ on an algebraic
surface $X$ is called {\it negative definite} if the intersection matrix
$(C_iC_j)$ of $C$ is negative definite.

(2)  A curve $C$ is called {\it rational} if for any effective
divisor $D=n_1C_1+\cdots+n_rC_r\neq 0$, we have $p_a(D)\leq 0$.

(3) If $A^2>0$, then the maximal reduced divisor
$C=C_1+\cdots+C_r$ with $CA=0$ is called the {\it exceptional
curve} of $A$. We denote it by $E(A)=C$.
\end{definition}

\begin{thm}\label{fundamental cycle} {\bf (Artin {\rm\cite{Artin1, Artin2}})}
Let $C=C_1+\cdots+C_r$ be a negative definite connected curve on
an algebraic surface $X$.

$(1)$ There is a unique effective divisor $Z=n_1C_1+\cdots+n_rC_r$
such that $ZC_i\leq 0$ and $Z$ is minimal. $(Z$ is called the {\bf
fundamental cycle} of $C).$ In fact, $Z\geq C$.

$(2)$ $p_a(Z)\geq 0$, and $p_a(Z)=0$ if and only if $C$ is
rational.
\end{thm}

\begin{thm} {\bf(Artin's projective contraction theorem {\rm\cite{Artin1, Artin2}})} A negative definite curve
$C=C_1+\cdots+C_r$ on a projective surface $X$ is rational if and
only if $C$ can be contracted to rational singular points on a
projective surface $Y$. $($The singular points on $Y$ is called rational if
$\chi(\mathcal O_X)=\chi(\mathcal O_{Y}))$.
\end{thm}

Benveniste \cite{Ben} generalized a result of Zariski
(\cite{Zar1}, Theorem 6.1).
\begin{thm}\label{Ben} {\bf (Benveniste {\rm\cite{Ben}})}  Suppose $C=C_1+\cdots+C_r$ is a
connected component of $E(A)$ for some nef and big divisor $A$. If
$C$ is rational and $n\gg0$, then $C$ can not be a fixed component
of $|nA|$.
\end{thm}

\section{Some technique results}

Reider's method is usually used to study the adjoint linear system
$|K_X+L|$ for a nef and big divisor $L$. In our case, $L$ is not necessarily
nef. Because there is no reference of this method for the general case,
we shall present in this section the generalization of Reider's method so that
Bogomolov's inequality can be used in the general case.

We use the notion ``$k$ points'' for any zero-dimensional
subscheme of length $k$, not requiring the points to be distinct.

Given a subscheme $Z'\subset Z$, the ``complement'' $Z''$ of $Z'$
in $Z$ is the canonical closed subscheme $ Z'' \subset Z$ with
an ideal sheaf $\I_{Z''}=[\I_{Z}:\I_{Z'}]$, i.e., for any open
set $U\subset X$,
$$
\I_{Z''}(U):=\{g\in \O_X(U) \ | \ g\I_{Z'}(U) \subset \I_Z(U)\}.
$$
We call $Z''$ the {\it residual subscheme} of $Z'$ in $Z$ and denote it by $$Z''=Z-Z'.$$

Assume that $Z$ is a local complete intersection, and $Z''$ is the
residual of $Z'\subset Z$ in $Z$. Then $Z'$ is the residual of
$Z''$ in $Z$. Furthermore, we have
$$\deg Z'+\deg Z''=\deg Z.$$

Note that in the surface case, the 4 equivalent conditions in the following theorem
imply that $\Delta$ is a local complete intersection.

\begin{theorem}\label{Theorem 2.5.2}  {Let $\Delta$ be a zero-dimensional
subscheme of $X$ $($including empty set$)$ and let $L$ be a divisor on $X$.
Then the following conditions are equivalent.
\begin{enumerate}
\item[{\rm(1)}]
     There is a rank two vector bundle $E$ with a non
zero global section $\delta$ satisfying
$$
Z(\delta)=\Delta, \hskip0.3cm \det E=L. \eqno(A)
$$
\item[{\rm(2)}] There are $3$ curves $F_1,F_2$ and $F_3$ such that
$F_1$ and $F_2$ have no common components, and
    $$
     \begin{cases}
       \Delta  =F_1\cap F_2-F_1\cap F_2 \cap F_3,    \cr
       L \equiv F_1+F_2-F_3.
     \end{cases} \eqno(B)
    $$
\item[{\rm(3)}]
There exists a rank two vector bundle $\mathcal E$ with a global
section $s$ such that $\dim Z(s)=0$ and
  $$
     \begin{cases}
       \Delta  =Z(s)-Z(s)\cap F,    \cr
       L \equiv \det\mathcal E-F.
     \end{cases} \eqno(C)
  $$
\item[{\rm(4)}]
     Either $\Delta=\emptyset$ or there is an element $\eta$ in
$H^{1}({\mathcal I}_\Delta(K_X+L))^\vee$
such that for any subscheme $($including empty set$)$
$\Delta'\subsetneq \Delta$,
$\eta$ is not in the image of the following natural inclusion map:
    $$
    H^{1}({\mathcal I}_{\Delta'}(K_X+L))^\vee \hookrightarrow
H^{1}({\mathcal I}_\Delta(K_X+L))^\vee.
    $$
Equivalently,
$$
\bigcup_{\Delta'\subsetneq \Delta}
H^1(\mathcal I_{\Delta'}(K_X+L))^\vee
\subsetneq
H^1(\mathcal I_{\Delta}(K_X+L))^\vee . \eqno(D)
$$
\end{enumerate}}
\end{theorem}
(See \cite{TanV} for the details of the proof).

\begin{remark}\label{Remark 2.5.3} In the above correspondence, if $\Delta=\emptyset$,
then the following trivial cases correspond to each other:
\begin{enumerate}
\item $E=\mathcal O_X\oplus \mathcal O_X(L)$;
\item $f=af_1+bf_2$ for some sections $a$ and $b$ of line bundles;
\item $f\in \text{im}(s)$ is in the image of $s$ (we do not go to
the details of this condition);
\item $\eta=0$.
\end{enumerate}
\end{remark}

We would like to mention the implication from (2) to (1) which will be used in the proof of
Lemma \ref{Lemma 2.7.9}.
Denote by $f_i$ the global section of $\mathcal O_X(F_i)$
defining $F_i$. Let $\mathcal F$ be the syzygy sheaf of
$(f_1,f_2,f_3)$,
\begin{equation}\label{equ3.1}
0\lra \mathcal F \lra \oplus_{i=1}^3\mathcal O_X(-F_i)
\overset{f}{\lra} \mathcal O_X,
\end{equation}
where $f$ is defined by $f(x,y,z)=f_1x+f_2y+f_3z$, and
let $E=\mathcal F(F_1+F_2)$. One can prove that $\det E=F_1+F_2-F_3=L $ and $E$ has a global
section $\delta$ such that
$$Z(\delta)=F_1\cap F_2-F_1\cap F_2\cap F_3=\Delta.$$

\begin{definition} We say that $\Delta$ satisfies Cayley-Bacharach property with
respect to $|K_X+L|$ if for any $F$ in $|K_X+L|$ and for any
subscheme $\Delta'\subset \Delta$ with $\deg\Delta'=\deg\Delta-1$,
$F$ contains $\Delta'$ implies that $F$ contains $\Delta$.
Equivalently, for any such $\Delta'$,
$$\dim H^1(\mathcal I_{\Delta'}(K_X+L))^\vee<\dim H^1(\mathcal I_\Delta(K_X+L))^\vee.
 \eqno(D')
$$
\end{definition}

$(D)$ implies $(D')$.

\begin{lemma} If $\Delta$ is reduced or $\deg\Delta\leq 2$, then $(D')$ is equivalent to
$(D)$.
\end{lemma}
\begin{proof} In the two cases, $\Delta$ admits at most a finite
number of subschemes $\Delta'$ with $\deg\Delta'=\deg\Delta-1$, so
$(D')$ implies $(D)$. Indeed, if $\Delta$ is reduced, the
finiteness is obvious. If $\Delta$ is a non-reduced
zero-dimensional scheme of degree 2, and if $p$ is a point on
$\Delta$, then it is easy to prove that $I_\Delta=(x,y^2)$, where
$x$ and $y$ are some local coordinates of $X$ near $p=(0,0)$. So
$\Delta$ contains only one subscheme $p$ of degree $1$.
\end{proof}

\begin{corollary}\label{Corollary  2.5.4} Let $L$ be a fixed divisor on $X$ and $k\geq 1$
a fixed integer. Then the following conditions are equivalent.
\begin{enumerate}
\item[{\rm(1)}]
$(A)$ \(equivalently $(B)$ or $(C)$\) has no solution for any
$\Delta\neq\emptyset$ with $\deg\Delta\leq k$.
\item[{\rm(2)}] For any zero dimensional subscheme $\Delta\neq\emptyset$ of
degree $\leq k$,
$$
H^1(\mathcal O_X(K_X+L))^\vee =H^1(\mathcal I_{\Delta}(K_X+L))^\vee.
$$
\item[\rm(3)] Any zero dimensional subscheme $\Delta\neq\emptyset$ of
degree $\leq k$ gives $\deg\Delta$ independent conditions on
$|K_X+L|$. Namely, $|K_X+L|$ is $(k-1)$-very ample.
\end{enumerate}
\end{corollary}

The author and E. Viehweg (\cite{Tan1, TanV, Tan2}) prove that
the Cayley-Bacharach theorem for an $n$-dimensional projective
manifold is equivalent to the $k$-very ampleness of $|K_X+L|$.

\begin{thm}\label{Theorem 2.3.8} {\bf (Bogomolov {\rm\cite{Bog}})}
Let $E$ be a rank two vector bundle on an algebraic surface $X$. If
$c_1(E)^2>4c_2(E)$, then there is an invertible subsheaf $\mathcal
O_X(M)\subset E$ such that
\begin{enumerate}
\item[{\rm(1)}] $(2M-c_1(E))H>0$ for any ample divisor $H$;
\item[{\rm(2)}] $(2M-c_1(E))^2\geq c_1^2(E)-4c_2(E)$;
\item[{\rm(3)}] for any nef divisor $A$,
   $$ MA\geq \dfrac12  c_1(E)A+\dfrac12\sqrt{A^2}\sqrt{c_1^2(E)-4c_2(E)} \ .$$
\end{enumerate}
\end{thm}
 (3) follows from (2) and Hodge index theorem.

\begin{lemma}\label{Lemma 2.6.1} Let $E$ be a rank two vector bundle on $X$, and
let $\M_1$ and $\M_2$ be two different maximal invertible
subsheaves of $E$. Then there exists an effective divisor $D$ on
$X$ such that
$$
c_1(E)-c_1(\M_1)-c_1(\M_2) \equiv D.
$$
Furthermore, $D=0$ if and only if $\F=\M_1\oplus\M_2$.
\end{lemma}

See \cite{Tan1} for the proof of this lemma.

\begin{theorem}\label{Theorem 2.6.2} Let $L$ be a divisor on $X$ such that $L^2>0$ and
$LH\geq 0$ for some ample divisor $H$. Assume that $\Delta$ is empty or a zero dimensional
subscheme of $X$ satisfying one of
the equivalent conditions of Theorem \ref{Theorem 2.5.2}. If
$L^2>4\deg\Delta$, then either
\begin{enumerate}
\item[{\rm(1)}] $\eta=0$ $($and so $\Delta=\emptyset$, this
corresponds to the trivial cases, see Remark \ref{Remark
2.5.3}$\,);$ or \item[{\rm(2)}] $\eta\neq0$, and there exists an
effective divisor $D\neq 0$ passing through $\Delta$ such that for
any nef and big divisor $A$,
$$
DL-\deg\Delta\leq D^2<\dfrac{\ell}{2}DA
\leq\dfrac{\ell}4\left(AL-\sqrt{A^2}\sqrt{L^2-4\deg
\Delta}\right),
$$
where $\ell=AL/A^2>0$.
\end{enumerate}
\end{theorem}
\begin{proof} From the assumption and Hodge index theorem,
we see that $LH>0$. By Riemann-Roch theorem, we can prove easily
that for a sufficiently large $n$, $h^0(nL)>0$. Hence for any nef
and big divisor $A$, $LA\geq 0$. Since $L^2>0$,  $LA>0$.

We assume that the equivalent condition (1) of Theorem
\ref{Theorem 2.5.2} is true. Namely there is a rank two vector
bundle $E$ with a non zero global section $\delta$ such that
$Z(\delta)=\Delta$, and $\det E=L$. Thus $E$ admits a maximal
invertible subsheaf $\mathcal M_1=\mathcal O_X$,
$$
0\longrightarrow \mathcal O_X\longrightarrow E\longrightarrow
\mathcal I_{\Delta}(L) \longrightarrow 0.
$$
From the assumption, $c_1^2(E)=L^2>4\deg\Delta=4c_2(E).$ Hence $E$
is not semistable. By Theorem \ref{Theorem 2.3.8}, $E$ admits
a new maximal invertible sheaf $\mathcal M_2=\mathcal O(M)$
satisfying the three inequalities in Theorem \ref{Theorem 2.3.8}. From
Theorem \ref{Theorem 2.3.8} (1) and $LH>0$ for
any ample divisor, we see that $M\neq 0$, so $\mathcal
M_1\neq\mathcal M_2$. Now by Lemma \ref{Lemma 2.6.1}, there exists
an effective divisor $D\equiv L-M $ passing through $\Delta$.
Substitute $M=L-D$ into the second and third inequality of
Theorem \ref{Theorem 2.3.8}, we get
\begin{align*}
DL&-\deg\Delta\leq D^2,\cr
DA&\leq\dfrac{1}2\left(LA-\sqrt{A^2}\sqrt{L^2-4\deg
\Delta}\right).
\end{align*}

By Lemma \ref{Lemma 2.6.1}, if $\eta\neq 0$, then $D\neq 0$. In
fact, we only need to prove that $D^2<\dfrac{\ell}{2} DA$. If
$DA=0$, then $D^2<0=\dfrac{\ell}{2} DA$ by Hodge index theorem. If
$DL>0$, also by Hodge index theorem and (2),
$$
D^2A^2\leq (DA)^2<DA\cdot LA/2,
$$
so $D^2<\ell DL/2$.
\end{proof}

\begin{corollary} {\bf (Beltrametti and Sommese {\rm \cite{BeS1})}}
As in Theorem \ref{Theorem 2.6.2}, if $L$ is nef and big and $\Delta\neq \emptyset$,
then we have
\begin{equation*}
DL-\deg\Delta\leq D^2< \dfrac12
DL<\deg\Delta,
\end{equation*}
\end{corollary}

\begin{lemma}\label{Lemma 2.7.2} Let $C=C_1+\cdots+C_r$ be a negative definite curve on $X$.
\begin{enumerate}
\item[{\rm(1)}] The classes of the $C_i$ are independent in $NS(X)\otimes \mathbb Q$.
\item[{\rm(2)}] Let $\mathcal E=x_1C_1+\cdots+x_rC_r$. If $\mathcal EC_i\leq 0$ $($resp.
$<0)$ for any $i$, then all $x_i$ are nonnegative $($resp.
positive$)$. $($Hence $\mathcal E$ is an effective divisor$)$.
\item[{\rm(3)}] If $A^2>0$, then the number of curves $C$ satisfying
$AC= 0$ is finite. Hence $E(A)$ is well-defined.
\item[{\rm(4)}] There is a nef and big divisor $A$ such that $E(A)=C$.
\end{enumerate}
\end{lemma}
\begin{proof}
(1) The proof is well-known.

(2) Write $\mathcal E$ in the form $\mathcal E=A-B$, where $A$ and
$B$ are effective divisors, without common components. We have
$\mathcal EB\leq 0$ by assumption, hence $AB-B^2\leq 0$. Since
$AB\geq 0$ and $B^2\leq 0$, it follows that $B^2=0$, and hence
$B=0$ since the subspace generated by $C_1,\cdots,C_r$ is
negative definite.

(3) By Hodge index theorem, these curves span a negative subspace
of $NS(X)\otimes \mathbb Q$. Thus the number is less or equal to
the dimension of $NS(X)\otimes \mathbb Q$.

(4) Let $H$ be a very ample divisor on $X$. Then we can find
integers $x_1,\cdots,x_r$ such that
$$|\det (C_iC_j)|\cdot HC_k+(x_1C_1+\cdots+x_rC_r)C_k=0,
\hskip0.3cm \text{ for } k=1,\cdots,r.
$$
By (2), $x_i$ are positive. Let $A=|\det (C_iC_j)| H +x_1C_1+
\cdots +x_rC_r$. Then $A$ is nef and big and $E(A)=C$.
\end{proof}

\section{Effective bounds}

In this section, we fix a nef and big divisor $A$ and an arbitrary
divisor $T$. Denote by $C_1,\cdots,C_r$ the exceptional curves of
$A$.

\begin{theorem}\label{Theorem 2.7.3} Let $A$ be a nef and big divisor, let
$T$ be any fixed divisor, and let $k$ be a nonnegative integer.
Assume that either $n>k+\frak M(A,T)$,  or $n\geq \frak M(A,T)$
when $k=0$ and $T\sim K_X+\lambda A$ for some $\lambda\in \mathbb
Q$. Suppose $|nA+T|$ is not $(k-1)$-very ample, i.e., there is a
zero dimensional subscheme $\Delta$ on $X$ with minimal degree
$\deg\Delta\leq k$ such that it does not give independent
conditions on $|nA+T|$. Then there is an effective divisor $D\neq
0$ containing $\Delta$ such that
\begin{equation}\label{equ4.1}
\begin{cases}
D^2+K_XD+k\geq TD, &\cr DA=0.
\end{cases}
\end{equation}
\end{theorem}
\begin{proof} Let $L=nA+T-K_X$, i.e., $|nA+T|=|K_X+L|$.
We claim that if $n>k+\frak M(A,T)$, or $k=0$, $n=\frak M(A,T)$
and $T\sim K_X+\lambda A$,  then
\begin{equation}\label{equ4.2}
L^2>4k , \hskip0.3cm
0\leq\dfrac12\left(AL-\sqrt{A^2}\sqrt{L^2-4k}\right)<1.
\end{equation}
Indeed,
$$
L^2-4k=A^2n^2+2A(T-K_X)n+(T-K_X)^2-4k,
$$
the bigger root of the above quadratic polynomial of $n$ is
$$
n_k=\left(-A(T-K_X)+\sqrt{(A(T-K_X))^2-(T-K_X)^2A^2+4kA^2}\right)/A^2.
$$
Let $h:=(A(T-K_X))^2-A^2(T-K_X)^2$. By Hodge index theorem, $h\geq
0$, with equality if and only if there is a rational number
$\lambda$ such that $T\sim K_X+\lambda A$. Since
$$
k+\frak
M(A,T)-n_k=\dfrac1{A^2}\cdot\left(1-\dfrac12\sqrt{h+4kA^2}\right)^2\geq0,
$$
we have $L^2>4k$ when $n>k+\frak M(A,T)\geq n_k$.

Let
$$f(x)=x^2-AL\cdot x+\dfrac{h}4+k\cdot A^2.$$
The smaller root of $f(x)$ is
$$
x_1=\dfrac12\left(AL-\sqrt{A^2}\sqrt{L^2-4k}\right).
$$
 On the other hand,
$$
f(0)\geq 0,\hskip0.3cm f(1)=A^2\cdot(k+\frak M(A,T)-n)< 0 \
\textrm{ when } n> k+\frak M(A,T),
$$
so $0\leq x_1<1$ when $n>k+\frak M(A,T)$.

Note that if $k=0$, $h=0$ and $n=\frak M(A,T)$, then $\frak
M(A,T)>n_0$ and $f(0)=f(1)=0$. So $x_1=0$ and (\ref{equ4.2}) is
also true. (In fact, this is the unique case where $x_1<1$ and
$f(1)\geq 0$.) This proves the claim.

 If $|nA+T|$ is not $(k-1)$-very ample, then there
exists a minimal zero dimensional subscheme $\Delta$ (may be
empty) with $\deg\Delta\leq k$ satisfying the conditions of
Theorem \ref{Theorem 2.5.2} corresponding to the non-trivial cases
(see Remark \ref{Remark 2.5.3}). Apply Theorem \ref{Theorem 2.6.2}
(2) to $\Delta$ and $L=nA+T-K_X$, we get an effective divisor
$D\neq 0$ containing $\Delta$ such that
\begin{align*}
&DL- \deg\Delta \leq D^2,\cr 0\leq DA
&\leq\dfrac12\left(AL-\sqrt{A^2}\sqrt{L^2-4\deg \Delta}\right)\leq
x_1 <1.
\end{align*}
It implies $DA=0$ and $D^2+K_XD+k\geq TD$.
\end{proof}

In general, if $x>0$ and $f(x)<0$, i.e.,
$$n>n_k+\dfrac{1}{A^2}\left(x+\dfrac{f(0)}{x}-2\sqrt{f(0)}\right)
=-\dfrac{A(T-K_X)}{A^2}+\dfrac{x}{A^2}+\dfrac{f(0)}{xA^2},
$$
then $x_1<x$, so $DA<x$.

The following are natural conditions on $A$ or $E(A)$ such that
(\ref{equ4.1}) has no nonzero solutions $D$:

\begin{enumerate}
\item[\rm(A)] $A$ is ample (see Corollary \ref{Matsusaka}).
(Matsusaka's condition, $\tau=+\infty$). \item[\rm(B)] $TC_i\geq
K_XC_i+k$, for any $i=1,\cdots,r$. (Laufer-Ramanujam's condition,
$\tau\geq k+1$). \item[\rm(C)] $E(A)$ is a rational curve.
(Artin's condition, $\tau(A,0)=2$).
\end{enumerate}

\begin{lemma}\label{Lemma 2.7.4} For any integer $k\geq 0$,
there are at most a finite number of effective divisors $D$
satisfying $(\ref{equ4.1})$.
\end{lemma}
\begin{proof} Let $C_1,\cdots, C_r$ be all of the curves
satisfying $AC_i=0$. We have known that these curves span a
negative definite subspace $W$ of $NS(X)\otimes \mathbb Q$. It is
easy to see that the first inequality in (\ref{equ4.1}) gives a bounded domain in
$W$. Thus if $D=\sum_{i=1}^rn_iC_i$ satisfies the conditions in
the lemma, then $(n_1,\cdots,n_r)$ must be in a bounded domain of
$\mathbb Q^n$. This implies the lemma since $n_i$ are integers.
\end{proof}

\begin{corollary}\label{Corollary 2.7.5} Assume that $A$ is nef and
big, and $k\geq 0$.
\begin{enumerate}
\item[{\rm(1)}]
If Laufer-Ramanujam condition is true for $k$ and $n>k+\frak
M(A,T)$, then $|nA+T|$ is $(k-1)$-very ample.
\item[{\rm(2)}]
If $E(A)$ is a rational curve, then $h^1(nA)=0$ for $n>\frak
M(A,0)$, and $|nA|$ is base point free for $n>1+\frak M(A,0)$.
\end{enumerate}
\end{corollary}
\begin{proof} This follows from Corollary 1.4.
\end{proof}

In general $T$ may not satisfy Laufer-Ramanujam condition. We will
modify it such that Laufer-Ramanujam condition is true.

 Let $\sigma_i=\max\{KC_i-T C_i+k,\ 0\}$ for
$i=1,\cdots,r$. Because $(C_iC_j)$ is a negative definite matrix,
we can find an integral divisor
$$
\mathcal E_k=\mathcal E_k(A, T, C_1+\cdots+C_r)= x_1C_1+\cdots+x_rC_r
$$
such that $\mathcal
E_kC_i=-|\det(C_iC_j)|\cdot\sigma_i$. By Lemma \ref{Lemma 2.7.2},
$\mathcal E_k$ is effective. Let $T'=T-\mathcal E_k$. Then
Laufer-Ramanujam condition for $k$ is true for $T'$, i.e.,
\begin{equation}
(T-\mathcal E_k)C_i\geq K_XC_i+k, \hskip0.3cm i=1,\cdots,r.
\end{equation}
By the previous corollary, when $n>k+\frak M(A,T-\mathcal E_k)$, $
|nA+T-\mathcal E_k| $ is $(k-1)$-very ample.

\begin{theorem}\label{Theorem 2.7.6} Assume that $A$ is nef and big. Let $T$ be
any divisor.
\begin{enumerate}
\item[{\rm(1)}] If $n>\frak M(A,T-\mathcal E_0)$, then
$$
h^1(nA+T)=h^1(\mathcal O_{\mathcal E_0}(nA+T)).
$$
Because $\mathcal O_{\mathcal E_0}(A)$ is a numerically trivial
bundle, by {\rm(\cite{CuS1}}, Theorem 8\,{\rm)},  $h^1(\mathcal O_{\mathcal
E_0}(nA+T))$ is a periodic function of $n$.
\item[{\rm(2)}] If $n>1+\frak M(A,T-\mathcal E_1)$, then the fixed part $B_n$
of $|nA+T|$ is bounded by $\mathcal E_1$,
and $B_n$ is a periodic divisor of $n$ by {\rm\cite{CuS2}}.
\item[{\rm(3)}] Let $C'$ be a connected component of $E(A)=C'+C''$, and let
$\mathcal E_1=\mathcal E_1'+\mathcal E_1''$ be the corresponding
decomposition $($here $T=0)$. If $C'$ is rational and $n>1+\frak
M(A,-\mathcal E_1'')$, then $C'$ can not be the fixed part of
$|nA|$. Namely, the fixed part of $|nA|$ is contained in $\mathcal
E_1''$. $($See Theorem \ref{Ben} or {\rm\cite{Ben}}$).$
\end{enumerate}
\end{theorem}
\begin{proof}
(1) Since $n>\frak M(A,T-\mathcal E_0)$, by the above corollary,
$|nA+T-\mathcal E_0|$ is $(-1)$-very ample, i.e.,
$H^1(nA+T-\mathcal E_0)=0$. In fact, if $n>\frak M(A,T-\mathcal
E_0)$, then $(K_X-nA-T+\mathcal E_0)A<0$, so
$$H^2(nA+T-\mathcal E_0)\cong H^0(K_X-nA-T+\mathcal E_0)^\vee=0.$$
From the long exact sequence of
the following
$$
0\longrightarrow \mathcal O(nA+T-\mathcal E_0)\longrightarrow
\mathcal O(nA+T)\longrightarrow \mathcal O_{\mathcal E_0}(nA+T)
\longrightarrow 0,
$$
we can see that
$$
h^1(nA+T)=h^1(\mathcal O_{\mathcal E_0}(nA+T)).
$$
Note that $\mathcal E_0$ is supported on these $C_i$ with
$C_iA=0$. It has been proved in \cite{CuS1} that $h^1(\mathcal
O_{\mathcal E_0}(nA+T))$ is a periodic function of $n$. This
completes the proof.

(2) We take $k=1$. If $n>1+\frak M(A,T-\mathcal E_1)$, then
$|nA+T-\mathcal E_1|$ is 0-very ample, so it has no fixed
component. Thus the fixed part $B_n$ of $|nA+T|$ is contained in
$\mathcal E_1$, i.e., $\mathcal E_1-B_n$ is effective. This
completes the proof.

(3) If $\mathcal E_1'=0$ or $\mathcal E_1''=0$, then (3) follows from (2) and
the previous corollary (2). Otherwise, we claim that $|nA-\mathcal E_1''|$ is 0-very ample.
Thus the fixed part of $|nA|$ is contained in $\mathcal E_1''$.

Indeed, with the notations of Theorem \ref{Theorem 2.7.3},
if we write $D=D'+D''$ and let $T=0$, then from
$$
(-\mathcal E_1)C_j\geq K_XC_j+1, \hskip0.3cm \textrm{ for all } j
$$
and $D'\cap D''=\emptyset$, we get
$$\begin{cases}
\mathcal E_1''D''+ K_XD''+1\leq 0, &\textrm{ if } D''\neq 0,\\
D'^2+K_XD'=2p_a(D'')-2\leq -2,  &\textrm{ if }  D'\neq 0,\\
D''^2\leq -1,  &\textrm{ if }  D''\neq 0.
\end{cases}$$
Since $D=D'+D''\neq 0$,
\begin{align*}
D^2+DK_X+1-(-\mathcal E_1'')D&=(D'^2+D'K_X)+D''^2+(D''K_X+1+\mathcal E_1''D'')<0.
\end{align*}
Now the claim is a consequence of Theorem \ref{Theorem 2.7.3}.
\end{proof}

\begin{theorem}\label{Theorem 2.7.7}
Assume that $A^2>0$ and $|A|$ has no fixed part. Let $T$ be any
divisor.
\begin{enumerate}
\item[{\rm(1)}] If $n>1+\frak M(A,0)$, then $|nA|$ has no base points.
$($Zariski {\rm\cite{Zar1}.)}
\item[{\rm(2)}] If $n>1+\frak M(A,0)$ and $n>1+\frak
M(A,T-\mathcal E_1)$, then $h^1(nA+T)=h^1(\mathcal O_{\mathcal
E_1}(T))$ is a constant, and $h^2(nA+T)=0$. So
$$h^0(nA+T)=\chi(nA+T)+s,$$
where $s=h^1(\mathcal O_{\mathcal E_1}(T))$ is a constant.
$($Castelnuovo {\rm\cite{Cas}).}
\item[{\rm(3)}] The fixed part $B_n$ of $|nA+T|$ is a fixed divisor for $n\gg 0$.
\end{enumerate}
\end{theorem}
\begin{proof} (1) Because $|A|$ has no fixed part, $A$ is nef. If $p$ is a
base point of $|nA|$, then there is a curve $D$ passing through
$p$ such that $DA=0$. Because $|nA|$ has also no fixed part, this
means that we can find a curve in $|nA|$ disjoint with $D$ because
$DA=0$. This is impossible since $p$ should be their common point.

(2) In this case, we can find a divisor $D$ in $|nA|$ such that
$D$ is disjoint with the exceptional curve $E(A)$ of $A$. Thus $\mathcal
O_{\mathcal E_1}(D)=\mathcal O_{\mathcal E_1}$. As in the proof of
Theorem \ref{Theorem 2.7.6},
$$
h^1(nA+ T )=h^1(\mathcal O_{\mathcal E_1}(nA+ T ))= h^1(\mathcal
O_{\mathcal E_1}(D) \otimes\mathcal O_{\mathcal E_1}( T
))=h^1(\mathcal O_{\mathcal E_1}( T )).
$$
Note that $n>1+\frak M(A,T-\mathcal E_1)$ implies $n>-A(T-\mathcal
E_1-K_X)/A^2 =-A(T-K_X)/A^2$, so $A(K_X-nA-T)<0$. Thus
$h^2(nA+T)=h^0(K_X-nA-T)=0$.

(3) Let $M_n$ be the moving part of $|nA+ T |$. Because $|A|$ has
no fixed part, $|M_n+A|$ has also no fixed part. Since $(n+1)A+ T
=M_n+A+B_n$, we have $B_{n+1}\leq B_n$. It implies that when
$n\gg0$, $B_n$ is a fixed divisor.
\end{proof}

Zariski's criterion (Theorem \ref{generation}) for the finite
generation of $R[D]$ gives a criterion for projective
contractability of a negative definite curve (see also
\cite{Sch}).

\begin{corollary}\label{Corollary 2.7.11} {\bf(Criterion for Projective Contractability)}
Let $C=C_1+\cdots+C_r$ be a negative definite curve. Then $C$ can
be contracted to {\rm(}normal singular{\rm)} points on a
projective surface if and only if there is a nef and big divisor
$A$ such that $C=E(A)$ and $|nA|$ has no fixed part for some $n$.
\end{corollary}

\begin{corollary}\label{Corollary 2.7.12} {\bf(Artin \cite{Artin1})}
A negative definite and rational curve $C$ on an algebraic surface
can be contracted projectively.
\end{corollary}
\begin{proof} By Lemma \ref{Lemma 2.7.2}, there is a nef and big divisor $A$ such that
$E(A)=C$. By Corollary \ref{Corollary 2.7.5}, $|nA|$ is base point
free for large $n$. Thus $C$ can be contracted projectively.
\end{proof}

\begin{theorem}\label{Theorem 2.7.8} Let $A$ be a nef and big divisor with exceptional curve
$E(A)=C_1+\cdots+C_r$, and let $T=0$. Denote by $E_1,\cdots,E_s$
the connected components of $E(A)$. Let $\mathcal E_{0,i}=\mathcal
E_0(A,0,E_i)$, let $Z_i$ be the fundamental cycle of $E_i$, and
let
$$
\widetilde{\mathcal E}_{0,i}=\begin{cases} \mathcal E_{0,i},
&\textrm{ if } \ \mathcal E_{0,i}\neq 0,\\
Z_i, &\textrm{ if } \ \mathcal E_{0,i}=0.\end{cases} \hskip1cm
\widetilde{\mathcal E}_0=\sum_{i=1}^s\widetilde{\mathcal E}_{0,i}.
$$
Assume that either $|A|$ has no base point, or $E(A)$ is rational.
If
$$n>2+\frak M(A,0),$$
 then
\begin{enumerate}
\item[{\rm(1)}] $\Phi_{nA}$ is a birational morphism onto a projective surface
$\Sigma_n=\Phi_{nA}(X)$.
\item[{\rm(2)}] On $X\setminus \cup_{i=1}^r C_i$, $\Phi_{nA}$ is an isomorphism.
\item[{\rm(3)}] $\Phi_{nA}$ contracts the curves $E(A)$ to some $($singular$)$
points of $\Sigma_n$.
\item[{\rm(4)}] Furthermore, if $E(A)$ is rational, then
$\Phi_{nA}$ has connected fibers. In general, if $n>\frak
M(A,-\widetilde{\mathcal E}_{0})$, then $\Phi_{nA}$ has also
connected fibers.
\end{enumerate}
\end{theorem}
\begin{proof} In fact, we only need to prove (4) that $\Phi_{nA}$ has connected
fibers. The locus over which $|nA|$ is not very ample is contained
in $E(A)$, so we only need to prove that
$\Phi_{nA}(E_i)\neq\Phi_{nA}(E_j)$ when $i\neq j$.

By construction, ${\mathcal E}_{0,i}=\mathcal E_0(A,0,E_i)=0$ if
and only if $E_i$ consists of $(-1)$- and $(-2)$-curves. By
definition, $\widetilde{\mathcal E}_{0,i}\neq 0$, and
$$
-\widetilde{\mathcal E}_{0}C_j\geq K_XC_j, \hskip0.3cm \textrm{
for } j=1,\cdots,r.
$$
So $|nA-\widetilde{\mathcal E}_0|$ satisfies Laufer-Ramanujian
condition for $k=0$. Hence it is $(-1)$-very ample, i.e.,
$H^1(nA-\widetilde{\mathcal E}_0)=0$.

Now we consider the exact sequence
\begin{equation}\label{short}
0\to\mathcal O(nA-\widetilde{\mathcal E}_0) \to \mathcal O(nA)\to
\mathcal O_{\widetilde{\mathcal E}_0}(nA)\to 0.
\end{equation}
Since $|nA|$ has no base point and $A\cdot E(A)=0$, the generic
divisor $B\in |nA|$ does not contain any $C_i$, and hence disjoint
with $E(A)$. We obtain
$$
\mathcal O_{\widetilde{\mathcal E}_0}(nA) = \mathcal
O_{\widetilde{\mathcal E}_0}(B)=\mathcal O_{\widetilde{\mathcal
E}_0}=\bigoplus_{i=1}^s\mathcal O_{\widetilde{\mathcal E}_{0,i}}.
$$
The long exact sequence of (\ref{short}) gives us a surjective map
$$
H^0(nA)\to \bigoplus_{i=1}^sH^0(\mathcal O_{\widetilde{\mathcal
E}_{0,i}}),
$$
so for each $i\neq j$, there is a section $s$ in $H^0(nA)$ such
that $s(\widetilde{\mathcal E}_{0,i})\neq 0$,
$s(\widetilde{\mathcal E}_{0,j})= 0$. Hence $|nA|$ separates $E_i$
and $E_j$, which means that $\Phi_{nA}(E_i)\neq \Phi_{nA}(E_j)$.

If $E(A)$ is rational, then one can prove similarly that
$H^1(nA-Z_i-Z_j)=0$ and $\Phi_{nA}(E_i)\neq \Phi_{nA}(E_j)$
provided
\begin{equation}\label{new}
n>\frak M(A,-Z_i-Z_j)=2+\frak M(A,0)-\dfrac{m_i+m_j}{4}
\hskip0.3cm \textrm{ for any } \ i\neq j,
\end{equation}
where $m_i=-Z_i^2$ is the multiplicity of the normal rational
singular point with exceptional curve $E_i$. On the other hand,
condition (\ref{new}) follows from our assumption $n>2+\frak
M(A,0)$.
\end{proof}

 Let $A$ be a divisor on $X$, and let
$$
R_m:=H^0(X,\mathcal O_X (mA)), \hskip0.3cm
R[A]:=\bigoplus_{m=0}^\infty R_m.
$$
$R[A]$ admits naturally a graded ring structure. The generation of
this ring was studied by Zariski \cite{Zar1} for any divisor $A$,
and by Mumford \cite{Mum}, Kodaira \cite{Kod1, Kod2} and Bombieri
\cite{Bom} for the canonical divisor of surfaces of general type.

\begin{lemma}\label{Lemma 2.7.9} Assume that $A^2>0$,
$\ell$ and $p$ are two positive integers such that $|\ell A|$ has
no base point, and $H^1(mA)=V$ is fixed for any $m\geq p$.
Let $k=\ell^2A^2$ and let
$$
\frak N(A,\ell,p):=\begin{cases} \max\left\{\,2\ell+p-1, \
3\ell+\dfrac{K_XA}{A^2} \,\right\}, & \text{ if \ } V=0, \cr
\max\left\{\, 2\ell+p-1,  \ k+\frak M(A,0)+\ell\, \right\}, &
\text{ if \ } V\neq 0,
\end{cases}
$$
If $m>\frak N(A,\ell,p)$, then we have
\begin{equation}\label{equ4.6}
R_m=R_\ell R_{m-\ell}.
\end{equation}
\end{lemma}
\begin{proof} Note first that $A$ is nef and big.
Let $E(A)=C_1+\cdots+C_r$. We choose three generic curves $F_1$,
$F_2$ and $F_3$ in $|\ell A|$ such that they have no common zero
point and do not contain any $C_i$. Then we see that $E(A)$ is
disjoint with all $F_i$. We denote by $f_i$ the global section
defining $F_i$. Let $\mathcal F$ be the syzygy sheaf of
$(f_1,f_2,f_3)$ (see (\ref{equ3.1})),
\begin{equation}\label{equ4.7}
0\lra \mathcal F \lra \oplus_{i=1}^3\mathcal O_X(-F_i) \lra
\mathcal O_X\lra 0,
\end{equation}
and let $E=\mathcal F(F_1+F_2)$. We have known that
$\det E=F_1+F_2-F_3=\ell A$ and $E$ has a global
section $\delta$ such that $Z(\delta)=F_1\cap F_2$,
\begin{equation}\label{equ4.8}
0\lra \mathcal O\lra E\lra \mathcal I_{F_1\cap F_2}(\ell A)\lra 0.
\end{equation}

From the conditions, we claim that
\begin{enumerate}
\item[a)] $H^1((m-\ell)A)=H^1((m-2\ell)A)=V$ if $m>2\ell+p-1$;
\item[b)] $ H^2((m-3\ell)A)=H^0(K_X-(m-3\ell)A)=0 $
if $m>3\ell+\dfrac{K_XA}{A^2}$;
\item[c)] $H^1(\mathcal I_{F_1\cap F_2}((m-\ell)A))=V$ if $V=0$ or
$m>\ell+k+\frak M(A,0)$.
\end{enumerate}
Indeed, we only need to prove c). From the long exact sequence of
$$
0\lra \mathcal O_X((m-3\ell)A)\lra \mathcal
O_X((m-2\ell)A)^{\oplus 2} \ \overset{(f_1,f_2)}\lra \ \mathcal
I_{F_1\cap F_2}(m-\ell) \lra 0,
$$
we obtain
$$
V^{\oplus 2}=H^1((m-2\ell)A)^{\oplus 2} \lra H^1(\mathcal
I_{F_1\cap F_2}((m-\ell)A)) \lra H^2((m-3\ell)A)=0.
$$
Thus if $V=0$, then $H^1(\mathcal I_{F_1\cap F_2}((m-\ell)A)=0$.
Now we consider the case $V\neq 0$. Suppose
$$
V=H^1(\mathcal O_X((m-\ell)A))\subsetneq H^1(\mathcal I_{F_1\cap
F_2}((m-\ell)A)).
$$
Then $F_1\cap F_2$ violates the $(k-1)$-very ampleness of
$(m-\ell)A$ ($k=\deg(F_1\cap F_2))$, so there exists a minimal non
empty subscheme $\Delta\subset F_1\cap F_2$ violating the
$(k-1)$-very ampleness. Because $m-\ell>k+\frak M(A,0)$, by
Theorem \ref{Theorem 2.7.3} we have a curve $D\neq 0$ passing
through $\Delta$ such that (\ref{equ4.1}) is true, hence the
support of $\Delta$ is contained in both $F_1$ and $E(A)$, which
contradicts the choice of $F_1$. This completes the proof of c).

Similarly, consider the long exact sequence of
(\ref{equ4.8})$\otimes \mathcal O((m-2\ell)A)$, we get
$$
V=H^1((m-2\ell)A)\to H^1(E((m-2\ell) A)) \to H^1(\mathcal I_{F_1\cap
F_2}((m-\ell)A))=V.
$$
So
\begin{equation}\label{II7-7}
h^1(\mathcal F(mA))=h^1(E((m-2\ell)A))\leq 2\dim V.
\end{equation}
The long exact sequence of (\ref{equ4.7})$\otimes\mathcal O_X(mA)$
gives us the following
\begin{align*}
&\oplus_{i=1}^3H^0((m-\ell)A) \ \overset{(f_1,f_2,f_3)}\lra \
H^0(mA) \lra H^1(\mathcal F(m)) \to \cr {\overset{\alpha}\to}\
&\oplus_{i=1}^3H^1((m-\ell)A) \lra \ H^1(mA) \lra H^2(\mathcal
F(m)).
\end{align*}
We can see by (\ref{II7-7}) that $\alpha$ is injective, equivalently
$(f_1,f_2,f_3)$ is surjective, namely
$$R_m=f_1R_{m-\ell}+f_2R_{m-\ell}+f_3R_{m-\ell}.$$
This completes the proof.
\end{proof}

The following theorem is about the projective normality of
$\Phi_{mA}(X)$ (see \cite{Bom}, Theorem 3A).

\begin{theorem}\label{Theorem 2.7.10} Let $A$ be a nef divisor with $A^2>0$ such that either
$|A|$ has no fixed part or $E(A)$ is rational. Let
 \begin{align*} \ell&=1+\frak
m(A,0),\\
p&=\begin{cases} \frak m(A,0), & \textrm{ if }
E(A) \text{ is rational, }\\
1+\max\{\frak m(A,0), \frak m(A,-\mathcal E_1)\}, & \textrm{ if }
|A| \textrm{ has no fixed part. }
\end{cases}
\end{align*}
Assume that
$$
2m>\begin{cases} \max\left\{\,2\ell+p+1, \
3\ell+3+\dfrac{K_XA}{A^2} \,\right\}, & \text{ if } E(A) \textrm{
is rational}, \cr \max\left\{\, 2\ell+p+1,  \ k+\frak
M(A,0)+\ell+1\, \right\}, & \text{ if } |A| \textrm{ has no fixed
part}.
\end{cases}
$$
Then for any $n\geq 1$,
$$
R_{nm}=R^n_m.
$$
\end{theorem}
\begin{proof} Under the conditions, $|\ell A|$ and $|(\ell+1) A|$
have no base point, and $h^1(mA)=V$ is fixed when $m\geq p$. By
the previous lemma, if $m>\frak N(A,\ell+1,p)$, then we have
$$
R_m=R_\ell R_{m-\ell}, \hskip0.3cm R_m=R_{\ell+1} R_{m-\ell-1}.
$$
Now we claim that for any $n\geq 1$,
$$
R_{nm}=R_m^n.
$$
Indeed, we can write $m=s\ell+t(\ell+1)$ for some non negative $s$
and $t$, (e.g., $s=\left[\frac{m}{\ell}\right](\ell+1)-m$ and
$t=m-\left[\frac{m}{\ell}\right]\ell$). We can assume that $n\geq
2$. From the assumption
$nm\geq 2m > \frak N(A,\ell+1,p)$, so the previous lemma
gives us that
$$
R_{nm}=R_{(n-1)m+s\ell+t(\ell+1)}=R_{(n-1)m}R_\ell^sR_{\ell+1}^t
=R_{(n-1)m}R_m.
$$
By induction on $n$, we have $R_{nm}=R_m^n$ for any $n\geq 1$.
\end{proof}

\begin{example} Let $\pi: X\to \mathbb P^2$ be a double cover
ramified over a smooth curve $B$ of degree $2d$, and let
$H=\pi^*(\mathcal O_{\mathbb P^2}(1))$ be the pullback divisor of
a line. Then $nH$ is very ample if and only if $n\geq d$. (Note
that if $n<d$, then $|nH|=\pi^*|\mathcal O_{\mathbb P^2}(n)|$, which implies that
the map defined by $|nH|$ factorizes through the
double cover $\pi$, so $nH$ is not very ample). On the other hand,
$K_X\equiv (d-3)H$ and $H^2=2$. Thus $\frak M(H,0)=d-\frac52,$
hence $\frak m(H,0)=d-2$. In particular, $n>2+\frak M(H,0)$ iff
$n\geq d$. Therefore, our bound in Corollary 1.3 can not be
improved for $X$ and $H$.
\end{example}

\end{document}